\documentclass{amsart}
\pdfoutput=1

\usepackage{amssymb}
\usepackage{microtype}
\usepackage{tikz}
\usetikzlibrary{arrows}
\usepackage{booktabs}
\usepackage[pdftex,colorlinks,citecolor=black,linkcolor=black,urlcolor=black,bookmarks=false]{hyperref}
\def\MR#1{\href{http://www.ams.org/mathscinet-getitem?mr=#1}{MR#1}}

\newtheorem{theorem}{Theorem}[section]
\newtheorem{proposition}[theorem]{Proposition}
\newtheorem{lemma}[theorem]{Lemma}

\newtheorem{fact}[theorem]{Fact}

\theoremstyle{definition}

\numberwithin{figure}{section}
\numberwithin{equation}{section}
\numberwithin{table}{section}

\newcommand{\R}{\mathbb{R}}
\newcommand{\C}{\mathbb{C}}

\newcommand{\E}{\mathbb{E}}
\newcommand{\prob}{\mathbb{P}}

\newcommand{\vol}{\operatorname{vol}}

\DeclareMathOperator{\tr}{tr}

\renewcommand{\Im}{\operatorname{Im}}

\title{Small-scale equidistribution for random spherical harmonics}

\author{Matthew de Courcy-Ireland}
\address{Department of Mathematics\\
Princeton University\\
Princeton NJ 08544} \email{mdc4@math.princeton.edu}

\date{October 26, 2017}

\thanks{We thank Peter Sarnak for his advice, encouragement, and support over the course of our work and NSERC for a PGS D grant.}

\begin{document}

\begin{abstract}
We study random spherical harmonics at shrinking scales. We compare the mass assigned to a small spherical cap with its area, and find the smallest possible scale at which, with high probability, the discrepancy between them is small simultaneously at every point on the sphere.
\end{abstract}

\maketitle

\section{Introduction} \label{sec:intro}

By \emph{random spherical harmonics}, we mean random functions $\phi : S^2 \rightarrow \R$ given by
\begin{equation*}
\phi = \sum c_j \phi_j
\end{equation*}
where the $2m+1$ functions $\phi_j$ form an orthonormal basis for degree $m$ spherical harmonics on $S^2$ and the coefficients are independent Gaussians of mean 0 and variance $1/(2m+1)$.
The choice of variance $1/(2m+1)$ guarantees that
if we integrate over a geodesic ball $B_r(z)$,
\begin{equation} \label{eqn:expvol}
\E\left[ \int_{B_r(z)} \phi^2 \right] = \frac{\text{vol}(B_r)}{4\pi} = \sin^2(r/2).
\end{equation}
In expectation, the random measure $\phi^2 d\text{vol}$ thus weights the ball $B_r(z)$ by its volume fraction. For an individual $\phi$, there is some deviation from the expected value, and this is our interest. All of our considerations apply equally well to complex-valued harmonics $\phi : S^2 \rightarrow \C$, replacing $\phi^2$ by $|\phi|^2$ throughout. Notice that the expected value in Equation~(\ref{eqn:expvol}) is independent of the center $z$, as it must be since the ensemble is invariant under rotation of $S^2$.

To normalize, consider the random variables
\begin{equation*}
X_z = \frac{1}{\text{vol}(B_r)} \int_{B_r(z)} \phi^2
\end{equation*}
so that
$\E\left[X_z\right] = \frac{1}{4\pi}$
is of order 1 
for all $r>0$ and $m \geq 1$. The \emph{discrepancy} is
\begin{equation*}
D(r,m) = \sup_{z} \left|X_z - \E[X_z] \right| = \sup_{z} \left| X_z - \frac{1}{4\pi} \right|. 
\end{equation*}

\begin{theorem} \label{thm:main}
If $r \rightarrow 0$ and $m \rightarrow \infty$ in such a way that
\begin{equation*}
\frac{rm}{\log{m}} \rightarrow \infty,
\end{equation*}
then for any fixed $\epsilon > 0$,
\begin{equation*}
\prob \{ D(r,m) > \epsilon \} \rightarrow 0.
\end{equation*}
\end{theorem}

In fact, the proof we give shows that
\begin{equation}
\prob \{ D(r,m) > \epsilon \} \leq C(\epsilon)m^2 e^{-c(\epsilon) rm }.
\end{equation}
for some positive constants $c(\epsilon)$ and $C(\epsilon)$, with $c(\epsilon)$ on the order of $\epsilon^2$. The hypothesis that $rm/\log{m} \rightarrow \infty$ guarantees that the factor $m^2$ can be absorbed, no matter how small a value $\epsilon$ is given. Thus the discrepancy $D(r,m)$ converges to 0 in probability as long as $rm \rightarrow \infty$ asymptotically faster than $\log{m}$. This means the random measure $\phi^2 d\text{vol}$ is approximately uniform at a scale $r \approx 1/m$, larger than $1/m$ by only a slowly growing function. The significance of $1/m$ is that it is the \emph{Planck scale}, namely $1/\sqrt{\lambda}$ where $\lambda = m(m+1)$ is the Laplace eigenvalue of any spherical harmonic of degree $m$. As the proof unfolds, we will see that $\phi^2 d\vol$ does not equidistribute below this scale. 
This is a quantum mechanical effect: There is enough mass but it is not being distributed evenly because the Planck scale sets a fundamental limit.\\

There is a heuristic justification of Theorem~\ref{thm:main} worth keeping in mind during the proof. To accurately sample a polynomial of degree $m$ requires a grid spacing of order $1/m$, and hence roughly $m^2$ points on $S^2$. With high probability, the maximum of $N$ independent Gaussians of unit variance is of order $\sqrt{\log{N}}$. Taking $N \asymp m^2$ and approximating the supremum by a maximum over $N$ points, we thus expect
\begin{equation*}
\sup_{z} | X_z - \E[X_z]] | = \sqrt{\text{var}} \sup_{z} \left| \frac{X_z - \E}{\sqrt{\text{var}}} \right| \approx \sqrt{\text{var}} \sqrt{\log{m}} .
\end{equation*}
In Section ~\ref{sec:variance}, we show that the variance is of order $1/(rm)$. So the discrepancy should be small when
\begin{equation*}
\frac{\log{m}}{rm} \rightarrow 0.
\end{equation*}

To make a rigorous proof out of the heuristic above, we need to be precise about approximating the supremum by a maximum over finitely many gridpoints. For a single point $z$, concentration of $X_z$ follows from the variance estimate in Section~\ref{sec:variance} provided only that $rm \rightarrow \infty$, no matter how slowly. To handle many points at once, we form a fine grid on the sphere $S^2$, and this is where it becomes necessary that $rm$ grow quickly enough. Suppose that the discrepancy satisfies $D(r,m) > \epsilon$, so that there is a point $z$ where $X_z$ deviates appreciably from its mean. If the grid is fine enough, then it is likely that a comparable deviation from the mean occurs also at a nearby gridpoint. In Section~\ref{sec:grid}, we gain control of how much $X_z$ and $X_{z_j}$ could differ when $z_j$ is a nearby gridpoint, and thus of how fine the grid must be. To control the probability that there is a deviant gridpoint, we simply use a union bound over the grid. Running this argument over a very fine grid loses a large factor, but we show in Section ~\ref{sec:chernoff} that the tails of $X_z$ at a single point are light enough to handle this loss provided that $rm \rightarrow \infty$ asymptotically faster than $\log{m}$.\\


In the course of proving the tail bounds, we give a fairly precise description of the random variable $X_z$. Expanding the square in $\phi^2$ shows that $X_z$ is a quadratic form in Gaussian random variables. This quadratic form can be diagonalized explicitly:
\begin{proposition} \label{prop:diag}
Fix a point $z \in S^2$ and consider $X = X_z$. If we choose for our basis functions $\phi_j$ the standard ultraspherical polynomials rotated so that $z$ is at the North pole $(0,0,1)$, then
\begin{equation*}
X = \sum \lambda_{\nu} \mathfrak{z}_{\nu}^2
\end{equation*}
where the $\mathfrak{z}_{\nu}$ are independent standard Gaussians for $0 \leq \nu \leq 2m$ and the coefficients $\lambda_{\nu}$ satisfy
\begin{equation}
\lambda_{2k} = \lambda_{2k+1} = \frac{1}{2 \pi^2} \sqrt{1 - (k/(rm))^2 }\frac{1}{rm}\left(1 + O_{\eta}\left(\frac{k^{2/3+\eta}}{rm}\right) \right).
\end{equation}
for any $\eta > 0$ and a ratio $0 \leq k/(rm) < 1$ bounded away from $1$. For $k \sim t$,
\begin{equation}
\lambda_k \ll_{\eta} (rm)^{-4/3 + \eta}.
\end{equation}
For $k$ so large that $k + k^p > t$, where $p > 1/3$, we have
\begin{equation}
\lambda_k \ll_p \frac{ \exp(-c k^{(3p-1)/2}) }{(rm)^2}
\end{equation}
for a constant $c > 0$.
\end{proposition}

\begin{figure}[h]
\includegraphics[width=0.5\textwidth]{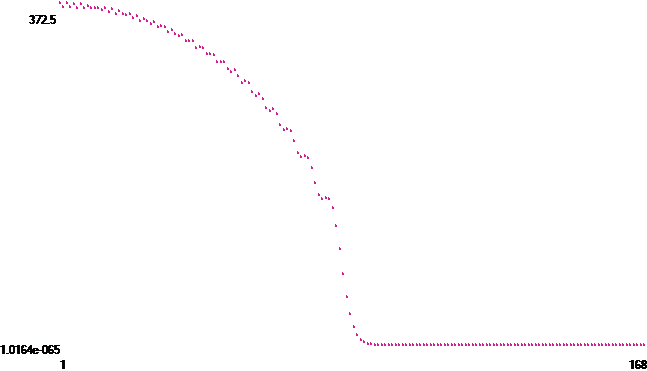}
\label{fig:besseldrop}
\caption{Dropoff of Bessel integrals $\int_0^{rm} uJ_k(u)^2du \div \int_0^{\pi m} uJ_k(u)^2 du$ for $0 \leq k \leq 2rm$ with $m=10000$, $r = \log(m)^2$. Made with pari-gp}
\end{figure}


As a consequence of Proposition \ref{prop:diag}, we control the tails of $X_z$:
\begin{lemma} \label{lem:tails}
For any $\epsilon > 0$ and any fixed $z \in S^2$, there are positive $C(\epsilon)>0$ and $c(\epsilon)>0$ such that
\begin{equation*}
\prob \left\{ \left| X_z - \frac{1}{4\pi} \right| > \epsilon \right\} \leq C(\epsilon) e^{-c(\epsilon)rm}.
\end{equation*}
The constant $c(\epsilon)$ in the exponent can be taken proportional to $\epsilon^2$.
\end{lemma}
This lemma gives exponential decay in $rm$, which is enough to absorb any power of $m$ sacrificed in tribute to the union bound because of the assumption that $rm$ is asymptotically larger than $\log{m}$.
We distinguish two cases in Lemma~\ref{lem:tails}: The upper tail where $X > \E[X] + \epsilon$ and the lower tail where $X < \E[X] - \epsilon$.  It seems natural to consider them separately because, at the lowest level of intuition, their origins are quite different. An easy way to imagine a lower tail event is that $B_r(z)$ contains a large part of the nodal set of $\phi$ so that the average of $\phi^2$ over $B_r(z)$ is small. On the other hand, the model scenario for an upper tail event is that $\phi$ achieves its maximum value near $z$ so that the average over $B_r(z)$ is large. Nevertheless, similar arguments in Section~\ref{sec:chernoff} control both the upper and lower tails.

Our original motivation for studying random spherical harmonics is the paper of Nazarov and Sodin on their nodal domains \cite{NS} and the far-reaching generalizations in \cite{NS2}. A natural context for Theorem~\ref{thm:main} is \emph{quantum unique ergodicity} (QUE). By QUE for a Riemannian manifold $M$, we mean that for any fixed measurable subset $A$ of $M$,
\begin{equation*}
\int_A \phi_{\lambda}^2 d\vol \rightarrow \vol(A)
\end{equation*}
for any sequence of Laplace eigenfunctions $\phi_{\lambda}$ with growing eigenvalue $\lambda \rightarrow \infty$. 
This is known to be false on $M = S^2$, because of the zonal spherical harmonics for example, but Rudnick and Sarnak conjecture that it is true on any compact negatively curved surface \cite{RS}. This has been shown for examples of arithmetic origin in work of Lindenstrauss \cite{L1}, \cite{L2}, and Bourgain-Lindenstrauss \cite{BouLi}, Jakobson \cite{J}, Holowinsky \cite{H}, Holowinsky-Soundararajan \cite{HS}. For progress constraining the possible limit measures in general, see Anantharaman \cite{A}, Anantharaman-Nonnenmacher \cite{AN}, Anantharaman-Silberman \cite{AS}, and Dyatlov-Jin \cite{DJ}. Even though QUE may fail for certain exceptional sequences of spherical harmonics, VanderKam \cite{VdK} shows that it does hold with probability tending to 1 for $\phi_{\lambda}$ in a randomly chosen orthonormal basis. Generating an entire basis at once is not the same as sampling from the monochromatic ensemble as we do here, but the two random models are similar. The scenario where the set $A$ shrinks as the eigenvalue grows has not been considered until recent papers such as Han-Tacy \cite{HT}, Granville-Wigman \cite{GW}, Lester-Rudnick \cite{LR}, Humphries \cite{Hum}.

\section{Some facts from analysis} \label{sec:besselfacts}

For ease of reference, here are some of the tools we use below.
\begin{fact} \label{fact:addition}
\textbf{(Addition formula for spherical harmonics)}
For any orthonormal basis of spherical harmonics $\phi_j$ of degree $m$, and for any points $x$ and $y$ on $S^2$,
\begin{equation}
\sum_j \phi_j(x)\phi_j(y) = \frac{2m+1}{4\pi} P_m(x \cdot y).
\end{equation}
Here, $P_m$ is the Legendre polynomial of degree $m$ normalized so that $P_m(1) = 1$. In particular, $| P_m | \leq 1$.
\end{fact}

\begin{fact} \label{fact:bernstein}
(\textbf{Bernstein's inequality})
The Legendre polynomial $P_m$ satisfies
\begin{equation}
P_m(\cos\theta)^2 \leq \frac{2}{\pi} \frac{1}{m \sin\theta}
\end{equation}
for all $\theta > 0$.
\end{fact}

\begin{fact} \label{fact:ultraspherical}
\textbf{(Basis of ultraspherical harmonics)}
Fix any point $z \in S^2$ as origin. There is an orthonormal basis of spherical harmonics of degree $m$ that are orthogonal not only over $S^2$ but also over any spherical cap $B_r(z)$ centered at $z$.
\end{fact}

In fact, the standard basis of ``$Y_l^m$"s has this property. Let the distance $\theta$ and the longitude $\alpha$ be spherical coordinates with respect to the point $z$. Then the $2m+1$ functions
\begin{equation} \label{eqn:basis}
\phi_{j,T} = \frac{ P_m^j(\cos\theta) T(j\alpha) }{ \left(\int_{0}^{2\pi} \int_0^{\pi} P_m^j(\cos\theta)^2 T(j\alpha)^2 \sin\theta d\theta d\alpha \right)^{1/2} } 
\end{equation}
form an orthonormal basis for spherical harmonics of degree $m$. The indices $j$ and $T$ run over $j = 0, 1, \ldots, m$ and $T \in \{ \sin, \cos \}$, excluding the case where $j=0$ and $T = \sin$, which gives $0$. These basis functions are orthogonal over any spherical cap $B_r(z)$ around $z$, no matter how small the radius $r$, because the functions $T(j\alpha)$ are orthogonal over the circle $0 \leq \alpha \leq 2\pi$. The polynomials $P_m^j$ are given by
\begin{equation*}
P_m^j(\cos\theta) = \frac{j!}{(2j)!} \frac{(m+j)!}{m!} (\sin\theta)^j P_{m-j}^{(j,j)}(\cos\theta).
\end{equation*}
in terms of Jacobi polynomials $P_n^{(\alpha,\beta)}$ with $n=m-j$ and $\alpha = \beta = j$.
We follow Szeg\H{o}'s treatment in section 4.7 of \cite{S}.
When $j = 0$, we have the Legendre polynomial of degree $m$. As $j$ increases, $P_m^j(x)$ vanishes to higher and higher order at $x=1$. This endpoint $x=1$ corresponds to the point $z$ on the sphere when we take $x = \cos{\theta}$, $\theta$ being the distance to $z$. 

\begin{fact} \label{fact:hilb}
(\textbf{Hilb asymptotics})
\begin{equation} \label{eqn:hilb}
P_m^j(\cos\theta) = h_{j,m} \left( \sqrt{\frac{\theta}{\sin\theta}} J_j((m+1/2)\theta) + O\left(\frac{(m-j)!}{m!} m^j \left(\frac{\sin\theta}{2}\right)^j\theta^{1/2} (m-k)^{-3/2} \right) \right)
\end{equation}
where
\begin{equation*}
h_{j,m} = \frac{j! 2^j}{(2j)!} \frac{(m+j)!}{(m-j)!} (m+1/2)^{-j}.
\end{equation*}
\end{fact}
The factor $h_{j,m}$ disappears when we normalize in $L^2$ and thus plays no role. Equation (\ref{eqn:hilb}) is a special case of Szeg\H{o}'s asymptotic (formula (8.21.17) in \cite{S}) for Jacobi polynomials $P_{n}^{(\alpha,\beta)}$. For $\alpha > -1$ and any real $\beta$, with $N = n + (\alpha + \beta + 1)/2$,
we have the estimate
\begin{equation*}
\left(\sin \frac{\theta}{2} \right)^{\alpha}\left(\cos \frac{\theta}{2} \right)^{\beta}P_{n}^{(\alpha,\beta)}(\cos \theta) = \frac{\Gamma(n+\alpha+1)}{n!} \sqrt{\frac{\theta}{\sin\theta}} \frac{J_{\alpha}(N\theta)}{N^{\alpha}} + \epsilon(n,\theta).
\end{equation*}
The error satisfies
\begin{equation*}
\epsilon(n,\theta) = \begin{cases}
                      \theta^{1/2}O(n^{-3/2}) \ \text{if} \ c/n \leq \theta \leq \pi^{-} < \pi \\
\theta^{\alpha+2}O(n^{\alpha}) \ \text{if} \ 0 < \theta \leq c/n
                     \end{cases}
\end{equation*}
for any fixed $\pi^-$ less than $\pi$ and any $c>0$, the implicit $O$ constants being subject to the choice of these parameters. In particular, $\epsilon(n,\theta) \lesssim \theta^{1/2}n^{-3/2}$ holds for all $\theta$.
In the special case where $\alpha = \beta = j$ and $n=m-j$, Szeg\H{o}'s asymptotic implies Fact~\ref{fact:hilb}. The case $\alpha = \beta = 0$ is Hilb's formula for Legendre polynomials, namely
\begin{equation*}
P_m(\cos\theta) = \sqrt{\frac{\theta}{\sin\theta}} J_0((m+1/2)\theta) + O\left(\frac{1}{m^{3/2}}\right)
\end{equation*}
For $k$ smaller than, say, $m/3$, we have $(1-k/m)^{-3/2} \leq 2$. For $k$ much smaller than $\sqrt{m}$, the factor $(m-k)! (m+1/2)^k/m!$ is also bounded. In that case, a consequence of equation (\ref{eqn:hilb}) is that (for $k$ much smaller than $\sqrt{m}$)
\begin{align*}
\frac{\int_0^r P_m^k(\cos\theta)^2 \sin\theta d\theta}{\int_0^{\pi} P_m^k(\cos\theta)^2 \sin\theta d\theta} &= \frac{\int_0^r \theta J_k((m+1/2)\theta)^2 d\theta + O\left( 2^{-k}(m-k)^{-3/2} r^k /k \right)}{\int_0^{\pi} \theta J_k((m+1/2)\theta)^2 d\theta + O\left( 2^{-k} (m-k)^{-3/2} k^{-1/2} \right)} \\
&= \frac{ \int_0^r \theta J_k((m+1/2)\theta)^2 d\theta }{\int_0^{\pi} \theta J_k((m+1/2)\theta)^2 d\theta } \left(1 + O\left( m^{-1/2}2^{-k} k^{-1/2}\right) \right) \\
&= \frac{ \int_0^{rm} xJ_k(x)^2 dx}{\int_0^{\pi m} xJ_k(x)^2 dx} \left(1 + O\left(m^{-1/2}2^{-k} \right) \right).
\end{align*}

Thus Hilb's formula naturally leads to the following integrals.

\begin{fact} \label{fact:bessel-integrals}
\textbf{(Some integrals involving Bessel functions)}
\begin{equation} \label{eqn:gr-5.54}
\int_0^t xJ_k(x)^2 dx = \frac{t^2}{2} \left(J_k(t)^2 - J_{k-1}(t)J_{k+1}(t) \right)
\end{equation}

\begin{equation} \label{eqn:basecase}
\int_{0}^{t} uJ_0(u)^2 du = \frac{1}{2}t^2 \left(J_0(t)^2 + J_1(t)^2\right).
\end{equation}
\end{fact}
This is formula 5.54 in \cite{GR}. It can be checked by differentiating both sides and using the recurrence relation between $J_k$, $J_k^{\prime}$, and $J_{k\pm 1}$. The second is formula (10.22.29) in the Digital Library of Mathematical Functions \cite{NIST}, and can be construed as the $k=0$ case of (\ref{eqn:gr-5.54}) with $J_{-1}=-J_1$. We don't use (\ref{eqn:gr-5.54}) in the proof, but we did use it to compute the integrals for Figure ~\ref{fig:besseldrop}. 

\begin{fact}\label{fact:bessel-asymptotics}
\textbf{(Asymptotics of $J$ Bessel functions)} 
For $k > x > 0$, we have
\begin{equation} \label{eqn:large-ord}
J_k(x) = (2\pi)^{-1/2} k^{-1/2} (1 - u^2)^{-1/4} e^{k (\sqrt{1 - u^2} - \sinh^{-1}(u^{-1}))}\left(1 + O\left(\frac{1}{\sqrt{x^2-k^2}}\right) \right)
\end{equation}
where $u = x/n$ is strictly between $0$ and $1$. For $x > k$, write $x = k\sec\beta$ with $0 < \beta < \pi/2$. Then
\begin{equation} \label{eqn:large-arg}
J_k(k\sec \beta) = \sqrt{\frac{2}{\pi k\tan\beta}}  \left( \cos(k(\tan\beta - \beta) - \pi/4) + O\left(\frac{1}{k\tan\beta}\right) \right)
\end{equation}
noting that $k\tan\beta = \sqrt{x^2 - k^2}$.
When $k$ and $x$ are too close, that is, $|x-k| < Ck^{1/3}$, these approximations become inaccurate and we use the upper bound
\begin{equation} \label{eqn:transition-bound}
J_k(x) \ll k^{-1/3}
\end{equation}
although it is possible to be much more precise.
\end{fact}
The first of these is formula 7.13.2 (14) in volume 2 of the Bateman Manuscript Project \cite{E}, page 87. Note that
\begin{equation*}
\frac{d}{du} \left( \sqrt{1-u^2} - \sinh^{-1}(u^{-1}) \right) = \frac{u^{-1}-u}{(1-u^2)^{1/2}} > 0
\end{equation*}
so the quantity in the exponent increases with $u$ from its limit $-\infty$ as $u\rightarrow 0$ to its value $-\log(1+\sqrt{2})$ at $u=1$.
The Bessel function $J_n(x)$ is exponentially small for small $x$ and oscillates with a decaying amplitude $\sqrt{2/(\pi x)}$ for large $x$. See formula 8.41(4) on p.244 of \cite{W} for equation (\ref{eqn:large-arg}). In between, there is a transition range of length $Cn^{1/3}$ centered at $x = n$. In this region, $J_n(x)$ achieves a maximum value of order $n^{-1/3}$ and also reaches its first positive zero. This maximum of order $n^{-1/3}$ is considerably larger than the amplitude $n^{-1/2}$ for $x$ beyond the transition range, and can be regarded as a ``boost" from the Airy function. The result, stated as 8.2(1) on p.231 of \cite{W}, is
\begin{equation*}
J_{n}(n) = \frac{\Gamma(1/3)}{2^{2/3} 3^{1/6} \pi} n^{-1/3} + O(n^{-2/3}).
\end{equation*}
In this regime, where $|x- n|$ is of order $n^{1/3}$ or smaller, Watson established an asymptotic for $J_{n}(x)$ stated as formulas (1) and (2) on p.249 of \cite{W} depending on which of $x$ and $n$ is the larger. Olver gives an asymptotic expansion for $J_{n}(n + \tau n^{1/3})$ in \cite{O}.


As a corollary of the behaviour of $J_{\nu}(t)$ for large $t$, we have
\begin{fact} \label{fact:squares}
\textbf{(Bessel version of $\sin^2 + \cos^2 = 1$)}
As $t \rightarrow \infty$,
\begin{equation*}
J_{\nu}(t)^2 + J_{\nu+1}(t)^2 \sim \frac{2}{\pi t}\left(1 + O_{\nu}\left(\frac{1}{t}\right)\right).
\end{equation*}
\end{fact}
We are imprecise about the dependence of the error term on $\nu$ because we only use it with $\nu = 0$ in connection with Equation~(\ref{eqn:basecase}).


\begin{fact} \label{fact:oscillatory}
If $f(y)$ is real-valued and continuously differentiable for $a < y < b$ with $f^{\prime}(y)$ positive and monotone, and $\inf f^{\prime} > 0$, then
\begin{equation}
\int_a^b e^{i f(y)} dy \lesssim \frac{1}{\inf f^{\prime} }
\end{equation}
\end{fact}
This is shown using integration by parts on p.124 of \cite{T}.


\section{An Exact Formula for the variance} \label{sec:variance}

\begin{lemma} \label{lem:variance}
For any point $z \in S^2$,
\begin{equation} \label{eqn:exact}
{\rm var}[X_z] = \frac{2}{\vol(S^2)^2} \int_{B_r(z)} \int_{B_r(z)} P_m(x \cdot x')^2 \frac{dx}{\vol(B_r)} \frac{dx'}{\vol(B_r)},
\end{equation}
where $P_m$ is the Legendre polynomial of degree $m$ normalized so that $P_m(1) = 1$. In particular,
\begin{equation} \label{eqn:estim}
{\rm var}[X_z]  \asymp \frac{1}{rm}.
\end{equation}
\end{lemma}
Equation (\ref{eqn:exact}) is an exact formula: It holds regardless of the relative sizes of $r$ and $m$. But if $rm \rightarrow \infty$, then (\ref{eqn:estim}) shows that the variance converges to 0. This is good enough for us to conclude using Chebyshev's inequality that at any point $z$
\begin{equation*}
\prob\left\{ |X_z - \E[X_z] | > \epsilon  \right\} \leq \frac{\text{var}(X_z)}{\epsilon^2} \lesssim \frac{1}{\epsilon^2} \frac{1}{rm} \rightarrow 0,
\end{equation*}
as long as $rm \rightarrow \infty$. For smaller $r$, the variance remains of order 1 or even diverges.\\


Before calculating the variance, let us verify that the mean is given by (\ref{eqn:expvol}): 
\begin{align*}
\E\left[ \int_{B_r(z)} \phi^2 \right] = \int_{B_r(z)} \E[ \phi^2] &= \int_{B_r(z)} \sum \phi_j(x)^2 \E[c_j^2] dx 
= \int_{B_r(z)} \frac{2m+1}{4\pi} \frac{1}{2m+1} 
\end{align*}
by linearity of expectation, expanding the square, and the fact that, for any orthonormal basis of harmonics $\phi_j$,
\begin{equation*}
\sum_{j} \phi_j(x)^2 = \frac{2m+1}{4\pi},
\end{equation*}
which follows from Fact~\ref{fact:addition}. Thus the expectation is the volume fraction, as claimed.

\begin{proof}
Let $\chi$ and $\chi^{'}$ be two functions on $M = S^2$. We imagine the indicator functions of two equal-sized balls $B$ and $B'$, but one could also use smooth cutoffs. The covariance we are interested in is
\begin{equation*}
{\rm Cov}\left[ \int_M \phi^2 \chi, \int_M \phi^2 \chi^{\prime} \right] = \E\left[ \int_M \phi^2 \chi \int_M \phi^2 \chi^{\prime} \right] - \E\left[ \int_M \phi^2 \chi \right] \E\left[ \int_M \phi^2 \chi^{\prime} \right].
\end{equation*}
With $\chi = \chi^{'}$, this becomes the variance. 

By linearity of expectation, and writing $dx$ instead of $d{\rm vol}(x)$,
\begin{align*}
&\E\left[ \int_M \phi^2 \chi \int_M \phi^2 \chi^{\prime} \right] = \int_M \int_M dxdx^{\prime} \chi(x)\chi^{\prime}(x^{\prime}) \E\left[ \left(\sum_j c_j \phi_j(x)\right)^2 \left(\sum_k c_k \phi_k(x^{\prime}) \right)^2 \right] \\
&= \int \int dx dx^{\prime} \chi(x) \chi^{\prime}(x^{\prime}) \sum_i \sum_j \sum_k \sum_l \phi_i(x)\phi_j(x)\phi_k(x^{\prime})\phi_l(x^{\prime}) \E[c_i c_j c_k c_l]
\end{align*}
Given four independent random variables $a, b, c,$ and $d$ with mean zero, the expectation $\E[abcd]$ is 0 unless the variables coincide, in which case we get the fourth moment $\E[a^4]$, or the variables are equal in pairs, in which case we get a product of variances $\E[a^2]\E[b^2]$ and so on. If all of the variables are Gaussian with mean zero and variance $\sigma^2$, the result is $3\sigma^4$ in the all-equal case or $\sigma^4$ in the equal-in-pairs case. So splitting the sum into the four cases \\
$\bullet i=j=k=l$, contributing $3\sigma^4 \phi_i(x)^2\phi_i(x')^2$\\
 $\bullet i=j \neq k =l$, contributing $\sigma^4 \phi_i(x)^2 \phi_k (x')^2$ \\
$\bullet i = k \neq j = l$ or $i = l \neq k = j$, each contributing $\sigma^4 \phi_i(x)\phi_i(x')\phi_j(x) \phi_j(x')$\\
shows that the quadruple sum is
\begin{equation*}
\sigma^4 \left( 3\sum_i \phi_i(x)^2\phi_i(x')^2 + \sum_{i\neq k} \phi_i(x)^2\phi_k(x')^2 + 2 \sum_{i\neq j} \phi_i(x)\phi_i(x')\phi_j(x)\phi_j(x') \right)
\end{equation*}
Since $3 = 1+2$, we can merge the first term into the second and third terms to provide the missing diagonal terms and factor the double sums into single sums:
\begin{equation*}
\left( \sum_i \phi_i(x)^2 \right)\left( \sum_k \phi_k(x')^2 \right) + 2 \left( \sum_i \phi_i(x) \phi_i(x^{\prime}) \right)\left( \sum_j \phi_j(x)\phi_j(x^{\prime})\right).
\end{equation*}
When we integrate, the first term will cancel with the product of expectations being subtracted in the definition of covariance. The second term can be expressed using the addition formula for spherical harmonics:
\begin{equation*}
\sum_{j} \phi_j(x)\phi_j(x') = \frac{2m+1}{4\pi} P_m(x \cdot x') = K_{\lambda}(x,x').
\end{equation*}
Here, $\lambda = m(m+1)$ is the eigenvalue of a degree $m$ harmonic for the spherical Laplacian on $S^2$.
The result is that
\begin{equation*}
{\rm Cov}\left[ \int_M \phi^2 \chi \ , \int_M \phi^2 \chi^{\prime} \right] = 2\sigma^{4} \int_M \int_M K_{\lambda}(x, x^{\prime})^2 \chi(x) \chi(x^{\prime}) dx dx^{\prime}.
\end{equation*}
For balls $B$ and $B^{\prime}$ in the sphere $S^2$, this becomes
\begin{equation*}
{\rm Cov}\left[ \int_B \phi^2 \ , \int_{B^{\prime}} \phi^2 \right] = \frac{2}{(4\pi)^2} \int_B \int_{B^{\prime}} P_m(x \cdot x^{\prime})^2 dxdx^{\prime}.
\end{equation*}
In particular, the variance is given by
\begin{equation*}
{\rm Var}\left[ \int_{B} \phi^2 \right] = \frac{2}{ {\rm vol}(S^2)^2 } \int_B \int_B P_m(x \cdot x^{\prime})^2 dx dx^{\prime}.
\end{equation*}
which establishes (\ref{eqn:exact}).
\end{proof}

There is another approach to proving the variance formula (\ref{eqn:exact}). The random variable $X = X_z$ is a quadratic form in Gaussians, so its moment generating function is explicit (see Equation (\ref{eqn:mgf})). Elsewhere, we use this to compute all moments of $X$ recursively and show that, when standardized to have mean 0 and variance 1, $X$ converges to a Gaussian as $rm \rightarrow \infty$. 
The higher moments are polynomials in the traces $\tr(A^p)$, where $A$ is the matrix with entries $\int_B \phi_j \phi_k d\vol/((2m+1)\vol(B))$. Using Fact \ref{fact:addition} repeatedly expresses this trace as a multiple integral of a product of Legendre polynomials, much like the second moment is expressed in terms of $P_m(x \cdot x')^2$. We have
\begin{equation}
\tr(A^p) = \frac{1}{(4\pi)^p} \int_B \cdots \int_B \prod_{i=1}^{p} P_m(x_i \cdot x_{i+1}) \frac{dx_1}{\vol(B)}\ldots \frac{dx_p}{\vol(B)},
\end{equation}
where the indices are taken cyclically so that $x_{p+1} = x_1$.

\begin{proof}
We turn to the proof of Equation~(\ref{eqn:estim}). We have Bernstein's inequality
\begin{equation*}
P_m(\cos{\theta})^2 \leq \frac{2}{\pi} \frac{1}{m \sin{\theta}} \leq \frac{1}{m\theta}
\end{equation*}
which improves on the trivial bound $P_m(\cos\theta)^2 \leq 1$ once $\theta = d(x,x') > 1/m$. Since $d(x,x')$ ranges all the way up to $2r$, if we assume that $rm \rightarrow \infty$, most values of $\theta$ appearing in the integral will enjoy a substantially improved bound on $P_m(\cos{\theta})$. Fix $x \in B_r(z)$. The points $x'$ lie in a ball $B_{2r}(x)$ around $x$, by the triangle inequality, and the integral of $P_m^2 \geq 0$ can only increase if we include all $x' \in B_{2r}(x)$ instead of only those in $B_r(z) \cap B_{2r}(x)$. Therefore, using spherical coordinates with respect to $x$ on $B_{2r}(x)$,
\begin{align*}
\int_{B_r(z)} \int_{B_r(z)} P_m(x \cdot x')^2 dx' dx &\leq \int_{B_r(z)} \int_0^{2\pi} \int_0^{2r} P_m(\cos\theta)^2 \sin\theta d\theta d\alpha dx \\
&\leq \int_{B_r(z)} 2\pi \int_0^{2r} \frac{2}{\pi} \frac{1}{m \sin\theta} \sin\theta d\theta dx \\
&= 8\frac{r \vol(B_r)}{m} \\
&\leq 2\pi \frac{\vol(B_r)^2}{rm}
\end{align*}
by Bernstein's inequality (Fact~\ref{fact:bernstein}). We also used $\vol(B_r) = 4\pi \sin(r/2)^2$ and $\sin(r/2) \geq r/\pi$ for $r \leq \pi$. Thus, by (\ref{eqn:exact}), $\text{var}[X_z] \leq C/(rm)$ with $C = 1/(4\pi)$.

The upper bound on $\text{var}[X_z]$ holds for any fixed $m$. To give a lower bound, we assume $rm \rightarrow \infty$. Then Hilb's asymptotics for $P_m$ show that this integral really is of order $(rm)^{-1} \text{vol}(B_r)^2$. 
Let $x \cdot x' = \cos\theta$, so $\theta = d(x,x')$, and let $\xi = d(z,x)$. By the triangle inequality, $B_{r-\xi}(x) \subset B_r(z)$. The integrand is nonnegative, so we have a lower bound
\begin{align*}
&\int_{B_r(z)} \int_{B_r(z)} P_m(x \cdot x')^2 dxdx' \\
&\geq \int_{B_r(z)} \int_{B_{r-\xi}(x)} P_m(x \cdot x')^2 dx' dx \\
&= \int_0^{2\pi} \int_0^r \int_{\alpha = 0}^{2\pi} \int_{\theta=0}^{r-\xi} P_m(\cos\theta)^2 \sin{\theta} d\theta d\alpha' \sin{\xi} d\xi d\alpha\\
&= (2\pi)^2 \int_0^r \int_0^{r - \xi} \left( \frac{\theta}{\sin\theta} J_0((m+1/2)\theta)^2 + O(m^{-3/2})\right) \sin\theta d\theta \sin\xi d\xi \\
&= (2\pi)^2 \int_0^r \int_0^{(r-\xi)(m+1/2)} uJ_0(u)^2 du (m+1/2)^{-2} \sin\xi d\xi + O\left(m^{-3/2}r^4\right) \\
&= (2\pi)^2 \int_0^r \frac{(r-\xi)^2 }{2} \left(J_0^2 + J_1^2 \right)((r-\xi)(m+1/2))\sin\xi d\xi + O(m^{-3/2} \vol(B_r)^2 )\\
\end{align*}
At this point, we restrict the range of integration further to $0 \leq \xi < (1-\delta)r$ so that $(r-\xi)(m+1/2)  \geq \delta rm$, which grows without bound by assumption. This allows us to use \ref{fact:squares}. The result is that
\begin{equation*}
 \int_{B_r(z)} \int_{B_r(z)} P_m(x \cdot x')^2 \frac{dx dx'}{\vol(B_r)^2} \geq \frac{(1-\delta)^2}{rm} \left(\frac{1}{2} - \frac{1-\delta}{3} \right) + O( (rm)^{-2} + m^{-3/2} )
\end{equation*}
Taking $\delta \rightarrow 0$, we have that for $rm \rightarrow \infty$,
\begin{equation}
\text{var}[X_z] \geq \frac{2}{(4\pi)^2} \frac{1}{6} \frac{1}{rm} \geq \frac{1}{480} \frac{1}{rm}
\end{equation}
We have used the Facts~\ref{fact:hilb},\ \ref{fact:bessel-integrals}, and \ref{fact:squares}.
\end{proof}
There is a factor of $12\pi$ between the crude upper and lower bounds above. One could use spherical trigonometry to evaluate the double integral more exactly, but upper and lower bounds of order $1/(rm)$ are all we need. We can also express the variance as $2\sum \lambda_j^2$ and use Proposition~\ref{prop:diag} to estimate the coefficients $\lambda_j$. See equation (\ref{eqn:variance-asymptotic}) below.

\section{Proof of Proposition~\ref{prop:diag}} \label{sec:ultra}

We fix $z \in S^2$ and use the basis from Fact~\ref{fact:ultraspherical}. The key advantage of this basis is that the off-diagonal entries of the matrix $A$ in $X = \mathfrak{z}^T A \mathfrak{z}$ all vanish. Thus
\begin{equation} \label{eqn:expansion}
X = \sum_{k=1}^{2m+1} \lambda_k \mathfrak{z}_k^2
\end{equation}
where each random variable $\mathfrak{z}_k$ is a standard Gaussian and there are no cross terms. The coefficients $\lambda_k$ are, for $1 \leq j \leq m$,
\begin{align*}
\lambda_1 &= \frac{1}{(2m+1)\vol(B_r)} \int_0^r P_m^0(\cos\theta)^2 \sin\theta d\theta \div \int_0^{\pi} P_m(\cos\theta)^2 \sin\theta d\theta \\
\lambda_{2j} &= \lambda_{2j+1} = \frac{1}{(2m+1)\vol(B_r)}  \int_0^r P_m^j(\cos\theta)^2 \sin\theta d\theta \div  \int_0^{\pi} P_m^j(\cos\theta)^2 \sin\theta d\theta.\\
\end{align*}
Our opening move is Hilb's formula:
\begin{align*}
\lambda_k &= \frac{1}{(2m+1)\vol(B_r)} \frac{\int_0^r P_m^k(\cos\theta)^2 \sin\theta d\theta}{\int_{0}^{\pi} P_m^k(\cos\theta)^2 \sin\theta d\theta} \\
&= \frac{1}{(2m+1)\vol(B_r)} \frac{ \int_0^{rm} xJ_k(x)^2 dx}{\int_0^{\pi m} xJ_k(x)^2 dx} \left(1 + O\left(\frac{1}{2^k m^{1/2}} \right) \right)
\end{align*}

To appraise the coefficients $\lambda_k$ with $k$ growing, we approximate the integral $\int_0^t xJ_k(x)^2 dx$ using Fact~\ref{fact:bessel-asymptotics}. Consider an initial range $x < k - k^p$, an intermediate range where $k - k^p< x < k + k^p$, and a final range where $k + k^p < x < t$. To begin, $0 < p < 1$. In the initial range $x < k$ so we change variables to $x = k \ \text{sech} \ \alpha$ and use equation (\ref{eqn:large-ord}). The lower limit $x=0$ corresponds to $\alpha \rightarrow \infty$ while the upper limit $x = k-k^p$ corresponds to $\alpha = \alpha_0 = \cosh^{-1}(k/(k-k^p)) \sim \sqrt{2} k^{(p-1)/2}$. This gives
\begin{equation} \label{eqn:initial-integral}
\int_0^{k- k^p} xJ_k(x)^2 dx \ll k \exp(2k(\tanh \alpha_0 - \alpha_0)) < \exp(-c k^{(3p-1)/2}),
\end{equation}
for some $c>0$, since $\tanh{\alpha} - \alpha \sim -\alpha^3 /3$ for small $\alpha$. The constant $c$ is positive and could be taken close to $2/3$. Thus (\ref{eqn:initial-integral}) shows that the initial range can be neglected as long as we choose $p > 1/3$.
Over the transition range, we have
\begin{equation} \label{transition-integral}
\int_{k - k^p}^{k + k^p} xJ_k(x)^2 dx \ll k^p k(k^{-1/3})^2 \ll k^{1/3 + p}.
\end{equation}
For large $x = k\sec\beta$, we have
\begin{align*}
xJ_k(x)^2 &= k\sec\beta \frac{2}{\pi k\tan\beta} \left( \cos^2(k(\tan\beta - \beta) - \pi/4) + O\left(\frac{1}{k\tan\beta}\right) \right)\\
&= \frac{1}{\pi \sin\beta} \left(1 + \sin(2k(\tan\beta - \beta) ) + O\left( \frac{1}{k\tan\beta} \right) \right). \\
\end{align*}
The change of measure $dx = k \sec\beta \tan\beta d\beta = d\beta k\sin\beta /\cos^2(\beta)$ cancels the $\sin\beta$ in the denominator above.
Thus on the final stretch of the integration,
\begin{equation*}
\int_{k+k^p}^t xJ_k(x)^2 dx =k  \int_{\sec^{-1}(1+k^{p-1})}^{\sec^{-1}(t/k)} \sec^2(\beta)\left(1 + \sin(2k(\tan\beta - \beta)) + O\left(\frac{1}{k\tan\beta}\right) \right) \frac{d\beta}{\pi} 
\end{equation*}
The lower limit of integration, $\sec^{-1}(1 + k^{p-1})$, is roughly 0. The $O$ term contributes $O(\log{k} + \log(t^2-k^2))$ when integrated by a change of variables $u = \tan\beta$:
\begin{align*}
\int_{\sec^{-1}(1 + k^{p-1})}^{\sec^{-1}(t/k)} \frac{1}{\tan\beta} \sec^2(\beta) d\beta &= \log \tan \beta \big]^{k=\sec^{-1}(t/k)}_{\sec^{-1}(1+k^{p-1})} \\
&= \log \sqrt{t^2/k^2 - 1} - \log \sqrt{ 2k^{p-1}+k^{2(p-1)} } \\
&= \frac{1}{2} \log(t^2 - k^2) - \frac{1}{2} (p \log{k} + \log{2} + \log(1 + k^{p-1}/2) ).
\end{align*}
This can be regarded as an error term as long as $t^2 - k^2$ is large. The term $\sec^2(\beta) \sin(2k(\tan\beta - \beta))$ oscillates enough to be of lower order when integrated. Indeed, change variables to $y = \tan \beta$, $dy = \sec^2(\beta) d\beta$ so that the integral is
\begin{equation*}
k \int_{\sqrt{2k^{p-1} + k^{2(p-1)}}}^{\sqrt{(t/k)^2-1}} \sin(2k(y - \arctan{y})) dy = k \Im \left[ \int_a^b e^{i f(y)} dy\right]
\end{equation*}
where $f(y) = 2k(y - \arctan{y})$, $b = \sqrt{(t/k)^2 - 1}$, and $a = \sqrt{2k^{p-1} + k^{2(p-1)}}$. We have
\begin{equation*}
f^{\prime}(y) = 2k \frac{y^2}{1+y^2}
\end{equation*}
which is positive and increasing, with a minimum value of $f^{\prime}(a) \asymp k^p$ on the interval of integration. It follows from Fact \ref{fact:oscillatory} that
\begin{equation*}
\int_a^b e^{i f(y) } dy \lesssim \frac{1}{f^{\prime}(a)} \lesssim k^{-p}
\end{equation*}
and therefore
\begin{equation*}
k \int_{\sqrt{2k^{p-1} + k^{2(p-1)}}}^{\sqrt{(t/k)^2-1}} \sin(2k(y - \arctan{y})) dy = O(k^{1-p}).
\end{equation*}

The main term is therefore
\begin{align*}
\frac{k}{\pi} \int \sec^2(\beta) d\beta &= \frac{k}{\pi} \tan\beta \big]^{\beta =\sec^{-1}(t/k)}_{\sec^{-1}(1+k^{p-1})} \\
&= \frac{1}{\pi} \sqrt{t^2 - k^2} + O(k\sqrt{(1+k^{p-1})^2-1}) \\
&= \frac{1}{\pi} \sqrt{t^2-k^2} + O(k^{(p+1)/2}) \\
\end{align*}
In order for this to be larger than our estimates for the initial range, we take $3p-1 > 0$. For the intermediate range to be smaller than the main term, we take $1/3 + p < 1$. Thus any exponent $1/3 < p < 2/3$ is allowed. For definiteness, we can take $p = 1/2$, although a value closer to $1/3$ would be more natural from the point of view of the transition for $J_k$.
Combining the three ranges shows that for $k < t$ (strictly, for $k+k^p < t$)
\begin{equation*}
\int_0^t xJ_k(x)^2 dx = \frac{1}{\pi} \sqrt{t^2 - k^2} + O\big(e^{-c k^{(3p-1)/2}} + k^{(p+1)/2} + \log{t} + k^{1-p} + \log(t^2 - k^2)\big).
\end{equation*}
We would like to take $p = 1/3$ to balance the powers of $k$, but the implicit constant diverges because of the initial range. However, we can choose $p$ slightly larger than $1/3$ to obtain, for any $\eta > 0$,
\begin{equation}
\int_0^t xJ_k(x)^2 dx = \frac{1}{\pi} \sqrt{t^2 - k^2} + O_{\eta}(k^{2/3 + \eta} + \log(t^2 - k^2)).
\end{equation}
When $k$ is slightly larger, so that $k - k^p > t$, only the initial segment contributes. In this case, the integral is dominated by $\exp(-c k^{(3p-1)/2})$ and is therefore negligible.
If $k - k^p < t < k + k^p$ so that the transition region contributes, the integral is still at most $O(k^{1/3 + p})$.

The coefficients at hand are given by a ratio of these integrals with $t = rm$ relative to $t = \pi m$. In the latter case, $t$ is always substantially larger than $k$ and we get
\begin{equation*}
\int_0^{\pi m} xJ_k(x)^2 dx = m + O_{\eta}(k^{2/3+\eta}).
\end{equation*}
The ratio is
\begin{equation*}
\frac{\int_0^{rm} xJ_k(x)^2 dx}{\int_0^{\pi m} xJ_k(x)^2 } = \frac{r}{\pi} \sqrt{1 - (k/(rm))^2} + O_{\eta}\left( \frac{k^{2/3 + \eta}}{m} \right).
\end{equation*}
When we incorporate the error from Hilb's formula, we get
\begin{equation}
\lambda_k = \frac{1}{(2m+1)\vol(B_r)} \left( \frac{r}{\pi} \sqrt{1 - \left(\frac{k}{rm}\right)^2 } + O_{\eta}\left(\frac{k^{2/3+\eta}}{m} + \frac{rm}{2^k m^{3/2}} \right) \right).
\end{equation}
That is, since each appears for two different basis functions ($\sin$ versus $\cos$)
\begin{equation}
\lambda_{2k} = \lambda_{2k+1} = \frac{1}{2 \pi^2} \sqrt{1 - (k/(rm))^2 }\frac{1}{rm}\left(1 + O_{\eta}\left(\frac{k^{2/3+\eta}}{rm}\right) \right).
\end{equation}
This explains the elliptical shape in Figure~\ref{fig:besseldrop}. Also, to leading order, the coefficients just for $k < rm$ are enough to match the expected value of $X$. Indeed,
\begin{equation}
\E\left[ \sum_{j < 2 rm} \mathfrak{z}_j^2 \lambda_j \right] \sim 2 \frac{1}{2\pi^2} \int_0^1 \sqrt{1 - u^2} du = \frac{1}{4\pi} = \E[X]
\end{equation}
up to an error of $O((rm)^{-1/3})$.
We also have, with $D$ the nearest integer to $rm$,
\begin{equation} \label{eqn:variance-asymptotic}
\text{var}\left[\sum_{k=1}^{D} \mathfrak{z}_k^2 \lambda_k \right] = \sum_k 2 \lambda_k^2 = \frac{\pi^{-4}}{D} \left( \int_0^1 1 - u^2 du + O(1/(rm)) \right) = \frac{2}{3\pi^4} \frac{1}{D} + O(D^{-2})
\end{equation}
which is another way to see that the variance is of order $1/(rm)$, as shown in Section \ref{sec:variance}, and even to find the constant of proportionality. Higher moments can likewise be expressed in terms of sums of powers of $\lambda_k$, and then estimated by integrals of $(1 - u^2)^{M/2}$.

\section{Union bound over a grid} \label{sec:grid}

Form a (deterministic) grid of points $z_j$ on $S^2$ such that every point is within $\delta$ of one of the gridpoints. If there is a point $z$ such that
\begin{equation*}
\left| X_z - \frac{1}{4\pi} \right| > \epsilon,
\end{equation*}
then we can expect the discrepancy to be high also for a nearby gridpoint. Indeed, if $d(z,z_j) < \delta$, then
\begin{equation*}
\epsilon < \left| X_z - \frac{1}{4\pi} \right| \leq \left| X_{z_j} - \frac{1}{4\pi} \right| + \left| X_z - X_{z_j} \right| \lesssim \left| X_{z_j} - \frac{1}{4\pi} \right|+ \frac{\delta}{r} \| \phi \|_{\infty}^{2}.
\end{equation*}
The last step follows from comparing integrals over two nearby balls as follows. 
For two sets $B$ and $B'$, we have
\begin{equation*}
\left| \int_B \phi^2 - \int_{B'} \phi^2 \right| \leq \int_{B \Delta B'} \phi^2 \leq \text{vol}(B \Delta B') \| \phi \|_{\infty}^2
\end{equation*}
For balls $B = B_r(z)$ and $B' = B_r(z')$, the volume of the symmetric difference depends both on $r$ and on the separation $\delta = d(z,z')$ between their centers. We have
\begin{equation*}
\text{vol}(B_r(z) \Delta B_r(z')) = O(\delta r)
\end{equation*}
by comparison with Euclidean rectangles, or by a more accurate calculation. 
Passing to averages, this gives
\begin{equation*}
\left|X_z - X_{z'} \right| \leq \frac{\text{vol}(B_r(z) \Delta B_r(z'))}{\text{vol}(B_r)} \|\phi \|_{\infty}^{2} \lesssim \frac{\delta}{r} \| \phi \|_{\infty}^{2}.
\end{equation*}
So either (writing $X_j$ for $X_{z_j}$) there is a $j$ such that
$|X_j - \frac{1}{4\pi}| > \epsilon/2$ \\
or $\| \phi \|_{\infty}^{2} \gtrsim r \epsilon/\delta$.\\
It follows from Th\'{e}or\`{e}me 7 in the paper \cite{BL} of Burq and Lebeau that $\| \phi \|_{\infty}$ is, with high probability, on the order of $\sqrt{\log{m}}$. Canzani and Hanin give another proof of this in \cite{CH}.
Thus the latter case where $\| \phi \|_{\infty}^2$ is at least of order $r\epsilon/\delta$ is very unlikely provided that we have a growing lower bound:
\begin{equation*}
\frac{\| \phi \|_{\infty}^{2}}{\log{m}} \gtrsim \frac{r \epsilon}{\delta \log{m}} \rightarrow \infty
\end{equation*}
as $rm \rightarrow \infty$. We can rewrite this in the form
\begin{equation*}
\frac{\| \phi \|_{\infty}^{2}}{\log{m}} \gtrsim \frac{rm}{\log{m}} \frac{\epsilon}{m\delta}
\end{equation*}
By hypothesis, $rm$ is asymptotically larger than $\log{m}$. So, for any fixed $\epsilon$, we can choose $\delta$ to be $1/m$. Then the probability of this case occurring will go to 0 as $rm \rightarrow \infty$.\\

For the former case, we have a union bound:
\begin{align*}
\prob \left\{ \exists j \ \left|X_j - \frac{1}{4\pi}\right| > \epsilon/2  \right\} &\leq \left(\text{number of points}\right) \prob \left\{ \left|X_1 - \frac{1}{4\pi} \right| > \epsilon/2 \right\} \\
&\lesssim \delta^{-2} \prob \left\{ \left|X_1 - \frac{1}{4\pi} \right| > \epsilon/2 \right\} \\
&\lesssim m^2 \prob \left\{ \left|X_1 - \frac{1}{4\pi} \right| > \epsilon/2 \right\}
\end{align*}
With $\delta = 1/m$ as above, we see that the union bound has cost us a factor of $m^2$, and we would pay an even steeper price of $m^d$ to apply it on a $d$-dimensional sphere instead of $S^2$. To afford it, we appeal to Lemma~\ref{lem:tails}. Since $rm/\log{m} \rightarrow \infty$, the bound $\exp(-c(\epsilon) rm)$ is $o(m^{-d})$ for any $d$.\\

In fact, Burq and Lebeau show that $\prob\{ \|\phi \|_{\infty} > c_0 \sqrt{\log{m}} + r \} \leq Ce^{-cr^2}$ for a specific constant $c_0$ and positive constants $C$ and $c$. In our context, this shows that the probability of the latter case is exponentially small with respect to $\epsilon^2 rm$. Thus it is no worse than the bound from Lemma~\ref{lem:tails} that we apply to the former case. The rate of convergence in Theorem~\ref{thm:main} is thus $O(\exp(-c(\epsilon) rm))$.


\section {Chernoff bound and proof of Lemma~\ref{lem:tails}} \label{sec:chernoff}

Lemma~\ref{lem:tails} is a special case of a more general fact about quadratic forms in Gaussians, which we state as
\begin{proposition} \label{prop:chernoff-gauss-quad}
If $\mathfrak{z}_j$ are independent Gaussians of mean $0$ and variance $1$ for $1 \leq j \leq D$, and the weights $\lambda_j \geq 0$ satisfy
\begin{equation}
\frac{A^-}{D} \leq \sum_{j=1}^{D} \lambda_j^2 \leq \frac{A^+}{D}
\end{equation}
and
\begin{equation} \label{eqn:lambda-max}
\max_{j=1}^{D} \lambda_j \leq \frac{M}{D}
\end{equation}
then the random variable $X = \sum_j \lambda_j \mathfrak{z}_j^2$ has exponential concentration as $D \rightarrow \infty$: For any fixed $\epsilon > 0$, there is a positive rate $c(\epsilon) > 0$ such that
\begin{equation}
\prob\{ | X - \E[X] | > \epsilon \} \leq \exp(-c(\epsilon) D).
\end{equation}
\end{proposition}

For example, if each $\lambda_j$ is $1/D$, then $X$ is a rescaled $\chi^2$ random variable with $D$ degrees of freedom, which exhibits concentration for large $D$. The role of the hypotheses is just to allow us to truncate the Taylor expansion of $\log(1 \pm x)$, and assumption (\ref{eqn:lambda-max}) could be relaxed to an upper bound on $\sum \lambda_j^3$.

\begin{proof} The Chernoff bound is
\begin{equation*}
\prob \{ X > \E[X] + \epsilon \} = \prob \left\{ e^{sX} > e^{s(\E[X]+\epsilon)} \right\} \leq \frac{\E\left[e^{sX}\right]}{e^{s(\E[X]+\epsilon)}}
\end{equation*}
where, given $\epsilon > 0$, the parameter $s$ is chosen to minimize the upper bound. Choosing $s = 0$ would give the trivial bound that probabilities are at most 1. Choosing an $s$ for which $\E[e^{sX}]$ is infinite would be even worse.
We write $X = \sum_j \lambda_j \mathfrak{z}_j^2$ for a quadratic form in Gaussian random variables $\mathfrak{z}_j$. In our case, the sum is indexed by $1\leq j \leq 2m+1$ and
\begin{equation*}
\lambda_j = \frac{1}{2m+1} \frac{1}{\vol(B_r)} \int_{B_r} \phi_j^2.
\end{equation*}
In general, we take $j \leq D$ as our indices and allow the coefficients $\lambda_j = \lambda_j(D)$ to depend on the number of variables.
The moment generating function can be computed explicitly. For $s \geq 0$ small enough that $1 - 2s\lambda_j > 0$ for all $j$,
\begin{equation} \label{eqn:mgf}
\E\left[e^{sX}\right] = \prod_{j} \left(1 - 2s\lambda_j \right)^{-1/2}
\end{equation}
since, by independence of the variables $\mathfrak{z}_j$, the quantity on the left factors as a product of Gaussian integrals. By differentiation, the optimal $s$ would solve
\begin{equation*}
\sum_j \frac{\lambda_j}{1 - 2s\lambda_j} = \E[X] + \epsilon.
\end{equation*}
Expanding the left in a geometric series gives
\begin{equation*}
\sum_{\nu=1}^{\infty} (2s)^{\nu-1} \sum_j \lambda_j^{\nu} = \E[X] + \epsilon.
\end{equation*}
Note that the first term $\nu = 1$ in the sum on the left is $\sum_j \lambda_j = \E[X]$, which cancels with the right. We may thus rewrite the equation for the optimal $s$ as
\begin{equation} \label{eqn:optimal-s-expansion}
\sum_{\nu=2}^{\infty} (2s)^{\nu-1} \sum_j \lambda_j^{\nu} = \epsilon.
\end{equation}
Any choice of $s$ gives some bound, and it is natural to choose $s$ by truncating this geometric series and solving the resulting equation. Keeping only the first term gives
\begin{equation*}
s_1 = \frac{\epsilon}{2} \frac{1}{\sum_j \lambda_j^2}.
\end{equation*}
One could keep two terms and solve a quadratic equation to get
\begin{equation*}
s_2 = \frac{\sum_j \lambda_j^2}{4\sum_j \lambda_j^3} \left( \sqrt{ 1 + 4\epsilon \frac{\sum \lambda_j^3}{\left( \sum \lambda_j^2 \right)} } - 1 \right)
\end{equation*}
which agrees with $s_1$ to first order in $\epsilon$. We will content ourselves with $s_1$. When we expand the logarithm, the terms of order $\epsilon^1$ cancel so that $s_1$ gives 
\begin{align*}
&\prob\{ X > \E[X] + \epsilon\} \leq \E[e^{s_1 X}] e^{-s_1 (\E[X] + \epsilon)} \\
&= \prod_j \left(1 - \frac{\epsilon}{\sum \lambda_j^2} \lambda_j \right)^{-1/2} \exp\left(-\epsilon \E[X]/ \left(2 \sum \lambda_j^2 \right)\right) \exp\left(-\epsilon^2/\left(2 \sum \lambda_j^2 \right)\right) \\
&= \exp\left(-\frac{\epsilon^2}{2\sum \lambda_j^2} - \epsilon \frac{\sum \lambda_j}{2\sum \lambda_j^2} - \frac{1}{2} \sum_j \log\left(1 - \epsilon \frac{\lambda_j}{\sum_k \lambda_k^2} \right) \right) \\
&= \exp\left(-\frac{\epsilon^2}{4\sum \lambda_j^2} + \sum_{\nu=3}^{\infty} \frac{1}{2\nu} \frac{\sum \lambda_j^{\nu}}{(\sum \lambda_j^2)^{\nu}} \epsilon^{\nu} \right)
\end{align*}
For $0 \leq x < 1/3$, we have the one-variable calculus exercise
\begin{equation*}
-\log(1-x) \leq x + \frac{3}{4}x^2.
\end{equation*}
Indeed, the claim follows for small $x$ from the series expansion for $\log$ and the range $x  < 1/3$ guarantees that the difference between the right and the left is in fact increasing.
If we can take $x = \epsilon \lambda_j / \sum \lambda_k^2$, which we will see shortly really is less than 1/3, then this will bound the product:
\begin{equation*}
\prod_j \left(1 - \frac{\epsilon \lambda_j}{\sum \lambda_k^2} \right)^{-1/2} \leq \exp\left( \frac{1}{2} \frac{\epsilon \E[X]}{\sum \lambda_k^2} + \frac{3}{8} \frac{\epsilon^2}{\sum \lambda_k^2} \right)
\end{equation*}
The terms that are first-order in $\epsilon$ cancel and the numbers have been rigged so that $3/8 - 1/2 = -1/8 < 0$, which gives a negative coefficient of $\epsilon^2$. The resulting bound is
\begin{equation*}
\prob \{ X > \E[X] + \epsilon \} \leq \exp\left( -\frac{\epsilon^2}{8} \frac{1}{\sum_j \lambda_j^2} \right)
\end{equation*}
Assuming that $\sum \lambda_j^2 \leq A_2 /D$, this implies that
\begin{equation*}
\prob \{ X > \E[X] + \epsilon \} \leq \exp(-c(\epsilon)D)
\end{equation*}
with $c(\epsilon) = \epsilon^2 /(8A_2)$ quadratic in $\epsilon$. Thus the probability of a deviation above the mean is exponentially small in $D$, as required. We claimed above that for each $j$, we may assume $\epsilon \lambda_j /\sum \lambda_k^2 < 1/3$ or, in other words, that $\lambda_{\max} < \frac{1}{3\epsilon} \sum \lambda_k^2$. One could certainly replace $1/3$ by any $\alpha < 1$ through a more vigorous Taylor expansion. The important point is that $\lambda_{\max}$ and $\sum \lambda_k^2$ have the same order of magnitude as $D \rightarrow \infty$, namely $1/D$. For if $\lambda_{\max} \leq M/D$ and $A^-/D \leq \sum \lambda_k^2$, then we will be guaranteed that $\lambda_{\max} < 1/(3\epsilon) \sum \lambda_k^2$ as long as $\epsilon < A^- /(3M)$ is sufficiently small (in absolute terms, with no reference to $D$).

For the lower tail, we rewrite $X < \E[X] - \epsilon$ as $-X > \E[-X] + \epsilon$ and apply the argument above to $Y = -X$. The details are slightly different because the moment generating function is now
\begin{equation}
\E[e^{sY}] = \prod_{k} (1 + 2s \lambda_k)^{-1/2}
\end{equation}
with a $1+2s\lambda_k$ instead of $1-2s\lambda_k$ in each factor. Thus any $s \geq 0$ is allowed and yields the bound
\begin{align*}
\prob \{ X < \E[X] - \epsilon\} &= \prob\{ Y > \E[Y] + \epsilon \} \\
&= \prob \{ e^{sY} > e^{s(\E[Y] + \epsilon)} \} \\
&\leq \E[e^{sY}] \exp(-s(\E[Y] + \epsilon)) \\
&= \exp(s(\E[X] - \epsilon)) \prod_k (1 + 2s\lambda_k)^{-1/2} \\
&= \exp\left( -\epsilon s + s\sum_k \lambda_k + \frac{1}{2}\sum_k - \log(1+2s\lambda_k) \right)
\end{align*}

The optimal $s$ would solve
\begin{equation*}
\sum_k \frac{\lambda_k}{1 + 2s\lambda_k} = \E[X] - \epsilon.
\end{equation*}
The first-order choice of $s$ is again
\begin{equation*}
s_1 = \frac{\epsilon}{2} \frac{1}{\sum_k \lambda_k^2}
\end{equation*}
although the second-order choice is different than in the case of the upper tail:
\begin{equation*}
s_2^{-} = \frac{\sum_k \lambda_k^2}{4 \sum_k \lambda_k^3} \left(1 - \sqrt{1 - 4\epsilon \frac{\sum_k \lambda_k^3}{\left( \sum_k \lambda_k^2 \right)^2}} \right).
\end{equation*}
Choosing $s_1$ and using the inequality $-\log(1+x) \leq -x + x^2/2$ for $x \geq 0$ gives an upper bound of
\begin{align*}
&\exp\left( -\epsilon s + s\sum_k \lambda_k + \frac{1}{2}\sum_k - \log(1+2s\lambda_k) \right) \\
&\leq \exp\left( -\frac{\epsilon^2}{2\sum_k \lambda_k^2} + \frac{\epsilon^2}{4(\sum_k \lambda_k^2)} \right) \\
&= \exp\left(-\frac{\epsilon^2}{4} \left(\sum_k \lambda_k^2\right)^{-1} \right) \\
&\leq \exp(-c(\epsilon) D).
\end{align*}
Thus the lower tail is also exponentially unlikely in $D$, provided only that $\sum \lambda_k^2 \leq A^+ D$.
\end{proof}

Another way to prove Propostion~\ref{prop:chernoff-gauss-quad} is to complexify and consider $\E[e^{it X}]$ instead of $\E[e^{sX}]$. Inverting the Fourier transform recovers the density of $X$. One can shift contours to show that the density is exponentially small away from $\E[X]$, but some care is needed in truncating the integral $\int_{-\infty}^{\infty} e^{-ix t} \E[e^{it X}] dt$ to a finite range $\int_{-T}^{T}$ and shifting the finite segment to an imaginary height $[-T,T]+iH$. The parameters $T$ and $H$ will both be small multiples of $D$, depending on the constants in the hypotheses, with sizes constrained relative to each other.\\

The sum $\sum \lambda_j^2$ is nothing but the variance that we saw in Section ~\ref{sec:variance}, which is of order $1/(rm)$. The largest coefficient $\lambda_{\max}$ is also of order $1/(rm)$, by integrating $P_m(\cos\theta)$ with the help of Hilb's formula as in Section ~\ref{sec:variance}. For $\epsilon$ small enough, we are thus guaranteed that $\epsilon \lambda_j / \sum \lambda_k^2 < 1/3$ for all $j$, as promised above. The argument above then applies, showing that the probability is exponentially small in $rm$. This is enough to overcome any factor $m^2$ or even a higher power coming from the union bound, as long as $rm$ is asymptotically larger than $\log{m}$.

\section{Conclusion} \label{sec:conc}

We have approximated the supremum
\begin{equation*}
\sup_{z \in S^2} \left| \frac{1}{\vol(B_r)} \int_{B_r(z)} \phi^2 - \frac{1}{4\pi} \right|
\end{equation*}
by a maximum over only finitely many points $z$. To control the error introduced this way, we made a brutish argument based on the union bound. 
We discuss a more sophisticated tool below, but the union bound is not as crude as it might seem. The exponentially light tail given by Lemma~\ref{lem:tails} is at the heart of why Theorem \ref{thm:main} is true. A helpful analogy is given by $k$ balls thrown at random into $n$ boxes, where one asks for the probability that each box receives close to $k/n$ balls as expected.

Dudley \cite{D} proved a general bound that applies to a separable, subgaussian process $X_t$ indexed by a metric space $(T,d)$. Normalizing so that $\E[X_t] = 0$ for convenience, the subgaussian assumption is that for all $\lambda \geq 0$,
\begin{equation*}
\E[e^{\lambda(X_s - X_t)} ] \leq e^{\lambda^2 d(s,t)^2/2}.
\end{equation*}
Dudley's conclusion is that
\begin{equation*}
\E\left[ \sup_{t \in T} X_t \right] \lesssim \int_0^{\infty} \sqrt{\log{N(T,d,\epsilon)}} d\epsilon,
\end{equation*}
where $N(T,d,\epsilon)$ is the smallest number of balls of radius $\epsilon$, in terms of the metric $d$, needed to cover $T$. The constant hidden inside $\lesssim$ is absolute and could be taken to be 12. This entropy method was used effectively by Feng and Zeldtich in \cite{FZ} and by Canzani and Hanin \cite{CH}. In applications, the metric $d$ is given by
\begin{equation*}
d(s,t) = \sqrt{ \E[(X_s- X_t)^2] },
\end{equation*}
and it is not quite a metric because it is possible to have $d(s,t)=0$ with $s \neq t$.
In our context of random spherical harmonics, $T = S^2$ is the sphere and
\begin{equation*}
X_z = X_z^{\pm} =\pm \left( \frac{1}{\vol(B_r)} \int_{B_r(z)} \phi^2 - \frac{1}{4\pi} \right).
\end{equation*}
The sign $\pm$ ensures that deviations above and below the mean can both be controlled.
Taking $\chi$ and $\chi'$ in the proof of Lemma~\ref{lem:variance} to be the indicator functions of the balls $B_r(z)$ and $B_r(z')$, we can express the (squared) metric $d(z,z')^2$ as
\begin{equation*}
\frac{4}{\vol(S^2)^4} \left( \int_{B_r} \int_{B_r} P_m(x \cdot x')^2 \frac{dx dx'}{\vol(B_r)^2} - \int_{B_r(z)} \int_{B_r(z')} P_m( x \cdot x')^2 \frac{dx dx'}{\vol(B_r)^2} \right).
\end{equation*}
By spherical symmetry, the first term $\int_{B_r}\int_{B_r}$ does not depend on the center of the ball $B_r$ while the second term $\int_{B_r(z)} \int_{B_r(z')}$ depends only on the spherical distance between $z$ and $z'$. We have $d(z,z')=0$, and indeed the first term exactly equals the second when $z = z'$. The first term is of order $1/(rm)$, as we saw in Lemma~\ref{lem:variance}. As $z$ and $z'$ become more distant, the second term decreases because of the decay of $P_m(x \cdot x')^2$ given, for example, by Fact \ref{fact:hilb}. It would be interesting to give another proof of Theorem~\ref{thm:main} by understanding the geometry of $S^2$ under this metric and, in particular, estimating the covering numbers $N(T,d,\epsilon)$.

We expect a proof using classical tools similar to those listed in Section~\ref{sec:besselfacts} to work for higher-dimensional spheres $S^d$ in place of $S^2$. We hope to prove an analogue of Theorem~\ref{thm:main} valid on any compact surface $M$. One can also ask for quantum limits on the bundle $S^{*}M$ instead of only the base manifold $M$, as in the full formulation of quantum unique ergodicity.




\end{document}